\documentclass[12pt]{article}
\usepackage{amsmath}
\usepackage{amsthm}
\usepackage{amsfonts}
\usepackage{eucal}
\usepackage{graphics}

\renewcommand{\baselinestretch}{1.2}
\newcommand{\field}[1]{\ensuremath{\mathbb{#1}}}
\newcommand{\Z}{\field{Z}}
\newcommand{\Q}{\field{Q}}

\newcommand{\calO}{\cal O}
\newcommand{\nD}{{\cal O}(nD)}

\DeclareMathOperator{\Pic}{Pic}
\DeclareMathOperator{\Div}{Div_{T}}
\DeclareMathOperator{\Hom}{Hom}

\newtheorem{theorem}{Theorem}[section]
\newtheorem{lemma}[theorem]{Lemma}

\theoremstyle{definition}

\theoremstyle{remark}

\newtheorem*{ack}{\bf Acknowledgement:}

\title{The ring of global sections of multiples of a line bundle 
on a toric variety\\{\small To appear in Proceedings of the AMS}}
\author{E. Javier Elizondo\thanks{Supported in part by grant CONACYT 3936-E}}
 
\begin{document}

\maketitle

\renewcommand{\baselinestretch}{-0.5}
\begin{abstract}
In this article we prove that for any complete toric variety, and for any
Cartier divisor, the ring 
of global sections of multiples of
the line bundle associated to the divisor is finitely generated.
\end{abstract}
\renewcommand{\baselinestretch}{1.5}

\section{Introduction}
Let us start by stating the result and a few consequences. Let
$X$ be an algebraic variety which is complete over a field $K$,
 and let $D$ be any effective Cartier
divisor on $X.$ We denote by $\nD$ the line bundle associated to $nD,$
with $n \geq 0.$
It is a very interesting problem to know if the ring 
\begin{equation}
\label{ring}
R\, := \,  \bigoplus_{n\geq 0} \, H^0 \,(X, \nD)
\end{equation}
is a finitely generated $K$-algebra. It was O. Zariski who worked out the case of algebraic surfaces in order to solve the Riemann-Roch problem, and gave examples where $R$ is not finitely generated, see \cite{zar-roch}.
In this  article we prove that the ring $R$ is finitely generated 
when $X$ is a complete toric variety (perhaps singular).   

The ring $R$ appears in many interesting problems, for example, 
if $X$ is a nonsingular projective variety, then
the finite generation of $R$ implies that the series 
\begin{equation}
\label{series}
\sum_{n \geq 0} \, \dim H^0 (X, \nD) \, t^n
\end{equation}
is rational. It was asked by S. D. Cutkosky and V. Srinivas  \cite{sri-zar} if this series is rational when $D$ is a nef divisor. Our result gives trivially a positive answer in the case of a complete toric variety, and $D$ any effective Cartier divisor.

Furthermore, 
the rationality of the series also allows us to compute, at least 
theoretically, the dimension
of  $ H^0 (X, \nD)$ in terms of $n$. This is the Riemann-Roch problem, see for example \cite{zar-roch} and \cite{sri-zar}. Observe that these dimensions 
are not given by
the Riemann-Roch theorem since the variety  $X$ can be singular, and the line bundle $\nD$ is not necessarily generated by its global sections.
We would like to make a final remark. The ring $R$ is a subring of the ring 
$S \, = \, \bigoplus_{D\geq{0}} \, H^0 (X, {\calO}(D))$, which is finitely generated in the case of toric varieties. In fact, D. Cox proves in \cite{cox-hom} that the ring $S$ is a 
polynomial ring graded by the the monoid of effective divisors classes in the Chow group $A_{k-1}(X)$ of $X$, where $k$ is the dimension of $X$.

\section{The ring of global sections}

In this section we prove that the ring $R$, which was defined by  equation (\ref{ring}),
is a finitely generated $K$-algebra. 
Throughout this section $X$ means a complete toric variety 
over a field $K,$ and $D$ 
a Cartier divisor in $X.$ The construction of the cone $C_R$, in the proof
of the Theorem, is a well know construction, see for example 
~\cite{vmh-batyrev}.\ 
We start by recalling an important lemma that
will be use in the proof of the theorem.
\begin{lemma}[Gordan]
\label{gordan}
If $\sigma$ is a strongly convex rational polyhedral cone, and \, ${\sigma}^{\vee}$ its dual,
then $S_{\sigma} \, = \, {\sigma}^{\vee} \, \cap \, M$ is a finitely generated semigroup.
\end{lemma}

Now, we are ready for the main result.
\begin{theorem}
\label{main}
Let $X$ be a complete toric variety, perhaps singular, and let $D$ be an 
Cartier divisor in $X$. Then the ring
$$
 R\, := \,  \bigoplus_{n\geq 0} \, H^0 (X, \nD)
$$
is  finitely generated  as a $K$-algebra.
\end{theorem}

\begin{proof}
We first observe that we can consider our divisor $D$ to be $T$-invariant.
For the following sequence is exact, see ~\cite[page 116]{dan-tova}
$$
0 \, \longrightarrow \, M \, \longrightarrow \, \Div \, X \, \longrightarrow \, \Pic \,(X) \longrightarrow \, 0
$$
where \, $\Div \, X$ \, is the group of $T$-invariant Cartier divisors and 
$M := {\Hom}_\Z (N, \Z)$ is the dual lattice of $N \cong {\Z}^k$. Then we can write  $D$ as \,$D = \, \sum_{i=1}^{s} \, a_i D_i$ with  $\{D_i\}$ the set of invariant divisors. We consider the convex rational polyhedron $P_{nD}$ in 
$M_{\Q} \, := \, M {\bigotimes}_{\Z} \Q$ defined as
$$
P_{nD} \, = \, \{u \, \in \, M_{\Q} \, | \, <u,v_i> \, \, \geq \, -n{a_i} \, 
\mbox{ for all } \, i \}
$$
where $v_i$ is the the first element in the lattice appearing in the 
divisor $D_i$. We know that generators for the space $H^0 (X, \nD )$ are given by the 
elements of\, $P_{nD} \, \cap \, M$, which is a finite set because $X$ 
is a complete variety, this also implies that $P_D$ is a  rational convex
polytope.  Let us embed the k-dimensional {\Q}\,-vector 
space ${M}_{\Q}$ 
into the 
(k+1)-dimensional {\Q}\,-vector space ${M}^{k+1}_{\Q} \, := \, {\Z}^{k+1} \bigotimes_{\Z} \, \Q$ as the hyperplane with equation $x_{k+1} \, = \, 1$, where $(x_1, \ldots , x_{k+1})$ are coordinates
for  ${M}^{k+1}_{\Q}$. Denote by $C_R$ the ($k+1$)-dimensional cone in ${M}^{k+1}_{\Q}$ generated by the rays starting at the origin and passing through the vertices of $P_D$ (see figure below).\vspace{1cm}

\hspace{4cm}
\includegraphics{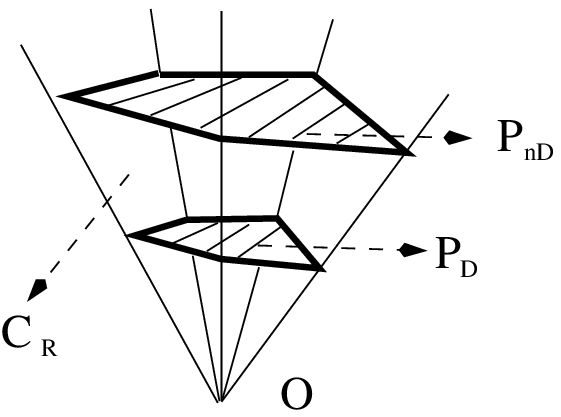}
\begin{center} 
{\em The cone \mbox{$C_R$}}
\end{center}

It follows 
from the definition that 
$P_{nD} = n P_{D}$, and this implies that the intersection of the hyperplane
$x_{k+1} = n$ with the cone $C_R$ is just the polyhedron $P_{nD}$. In other words,
the cone $C_R$ is the cone associated to the ring $R$, in the sense that any 
element of the ring $R$ is a finite linear combination of 
  integral points of $C_R \cap {\Z}^{k+1}.$
 The theorem follows  since the semigroup  $C_R \cap {\Z}^{k+1}$ is
finitely generated by Lemma ~\ref{gordan} (Gordan).
\end{proof}

\begin{ack}
I thank  very much V. Srinivas for suggesting me this problem  and for many 
wonderful and
fruitful conversations. I also would like to thank A. King for bringing 
to my attention the article of D. Cox, and to the referee for 
pointing out the article of V. Batyrev. \end{ack}



\begin{minipage}[t]{3in}
{ Instituto de Matem\'aticas \vspace{-0.25cm}\newline
Ciudad Universitaria, UNAM \vspace{-0.25cm}\newline
M\'exico D.F. 04510} \vspace{-0.25cm}\newline
M\'exico
\end{minipage}\vspace{.3cm}

{  e-mail: javier@math.unam.mx}

\end{document}